\documentclass{amsart}
\usepackage{ae,lmodern}
\usepackage[utf8]{inputenc}
\usepackage[T1]{fontenc}
\usepackage{graphicx}
\usepackage{amscd}

\def\stb{\underset{stb}{\sim}}
\def\01stb{\underset{0,1}{\sim}}

\def\f{\varphi}
\def\r{\rangle}
\def\l{\langle}
\def\Z{\mathbb{Z}}

\def\H{\mathbb{H}}

\def\sn{\mathrm{sn}}

\begin{document}

\newtheorem{pro}{Proposition}[section]
\newtheorem{theo}[pro]{Theorem}
\newtheorem{corol}[pro]{Corollary}
\newtheorem{lem}[pro]{Lemma }
\newtheorem{defi}[pro]{Definition}
\newtheorem{prob}[pro]{Problem}
\newtheorem{conj}[pro]{Conjecture}
\newtheorem{quest}[pro]{Question}
\newtheorem{rem}[pro]{Remark}
\newtheorem{remarques}[pro]{Remarks}
\newtheorem{exs}[pro]{Exemples}

\title{Some criteria for stably birational equivalence of quadratic forms}%
\address{INSPE de Franche-Comte, Site de Belfort, 
49 faubourg des anc\^etres, 90 000 Belfort}%
\email{sylvain.roussey@univ-fcomte.fr}%
\author{Sylvain Roussey}


\date{\today}%
%\dedicatory{To the ones I love}%
%\commby{}%
% ----------------------------------------------------------------
\begin{abstract} Let $\f$ and $\psi$ be quadratic forms over a field $K$ of characteristic different from 2. In this paper, we  give a criterion for isotropy of $\f$ over the function field of $\psi$ in terms of representations and we apply it to stably birational equivalence of $\f$ and $\psi$.
Then we use this criterion to investigate the case of stably birational equivalence of multiples of Pfister forms. Finally, we give a a criterion of stably birational equivalence of quadratic forms in terms of isomorphisms of quotients of special Clifford groups by the kernel of the spinor  norm.
\end{abstract}
\maketitle
% ----------------------------------------------------------------

\section{Introduction}

Let $\f$  be
a nondegenerate quadratic form over a field $K$ of characteristic different from 2. 
The function field of $\f$, denoted by $K(\f)$, is the function field of the  
quadric defined by the equation $\f=0$. 
An important and difficult question in the algebraic
theory of quadratic forms asked by  Lam in \cite{lam2}, 
chapter X, question 4.4B is the following  : 

\begin{prob}\label{iso}Consider a
quadratic form $ \f $. For which $K$-quadratic forms
$ \psi $ is $ \f_{K (\psi)} $ isotropic?\end{prob}

We first present a
characterization of isotropy analogous to Knebusch's norm theorem and
to Pfister's third representation theorem. This criterion is equivalent
to another characterization due to Witt (\cite{witt}, 1956), Knebusch, (\cite{norme}, 1973) and Bayer-Fluckiger (\cite{eva},
1995). More precisely, if we denote by $ T_{K} (\f) $ the multiplicative group generated by $D_{K}(\f)$, the set of values represented by $\f$ and if we denote by $ N_{K} (\f) $ the  multiplicative group generated 
 by $ D_{K} (\f) D_{K} (\f) $ in
$ T_{K} (\f) $, we get the following result :   \vspace*{0.2cm}

\noindent \textbf{Theorem  \ref{iocentrik}.} \itshape  Let $ \f $ and $ \psi $ be two anisotropic $ K $- quadratic forms.\\
If we write 
$ X = X_{1}, \ldots, X _{\dim (\psi)} $, then the following assertions
 are equivalent : 
\begin{enumerate}
\item $ \f_{K (\psi)} $ is isotropic ; 
\item $ D_{F} (\psi) D_{F} (\psi) \subset D_{F} (\f) D_{F} (\f) $ for 
every field extension $F/K$  ; 
\item $ a \psi (X) \in
D_{K (X)} (\f) D_{K (X)} (\f) $ for
every $ a \in D_{K} (\psi) $ ; 
\item $ N_{F} (\psi) \subset N_{F} (\f) $ for 
every field extension $F/K$ ; 
\item $ a \psi (X) \in N_{K (X)} (\f) $
for every $ a \in D_{K} (\psi) $ ; 
\item $ T_{F} (a \psi) \subset T_{F} (\f) $ for 
every field extension $F/K$ and every $ a \in D_{K} (\psi) $ ; 
\item $ a \psi (X) \in T_{K (X)} (\f) $
for every $ a \in D_{K} (\psi) $.
\end{enumerate}
\normalfont \vspace*{0.2cm}

More generally, considering the function field of an affine hypersurface, we get the following result in collaboration with Detlev W. Hoffmann : \vspace*{0.2cm}

\noindent \textbf{Theorem  \ref{coroeva}.} \itshape 
Let $ \f $ be an anisotropic quadratic form of dimension $n$ which
represents 1. We write $X=X_1,\ldots,X_n$. Let $ f \in K [X] $ be an irreducible polynomial
 where all $ X_{i} $ appear. If the degree of the term of highest degree for the lexicographical order
 is $ M $ and if $ f $ is unitary for the lexicographical order, then the following assertions are equivalent
 :  \begin{enumerate} \item  $ f \in D_{K(X)(\f)^m}$ with $m\leq M$ ; 
\item $ \f_{K (f)} $ is isotropic. \end{enumerate}\normalfont \vspace*{0.2cm}

Moreover, since  problem  \ref{iso} concerns an order relation, it
 is natural to wonder about the induced equivalence relation. So we say that $ \f $ and $ \psi $ are stably
birationally equivalent and we write $ \f \stb \psi $ if both
$ \f_{K (\psi)} $ and $ \psi_{K (\f)} $ are isotropic, and a new
question appears.  

\begin{prob}\label{stb} For which
$ K $-quadratic forms $ \f $ and $ \psi $ do we have $ \f \stb \psi $?
\end{prob}

Theorem \ref{iocentrik} easily implies the following characterisation of stably birational equivalence of quadratic forms  :   \vspace*{0.2cm}

\noindent \textbf{Theorem  \ref{zoumbalawe}.} \itshape 
Let $ \f $ and $ \psi $ be two anisotropic $K$-quadratic forms of dimension $\geq2$. The following assertions are equivalent : 
\begin{enumerate}
\item $ \f \stb \psi $ ; 

\item for 
every field extension $F/K$, $ D_{F} (\psi) D_{F} (\psi) = D_{F} (\f) D_{F} (\f) $ ; 

\item for 
every field extension $F/K$, $ N_{F} (\psi) = N_{F} (\f) $.
\end{enumerate}

\noindent If  $ \f $ and $ \psi $ represent 1, these assertions
are also equivalent to
\begin{enumerate}
\item [(4)] for 
every field extension $F/K$, $ T_{F} (\psi) = T_{F} (\f) $. \end{enumerate}
\normalfont \vspace*{0.2cm}

In \cite{unbelareticledefotzegraldquejeneconnaisaispas}, theorem
3.2, Fitzgerald has investigated the link between hyperbolicity of $\f$ over $K(\psi)$ and 
the hyperbolicity of $\pi\otimes\f$ over $K(\pi\otimes\psi)$ where $\pi$ is a $n$-fold  Pfister form. Elman and Lam gave implicitly the same kind of criterion for isotropy (\cite{u1}, proposition 2.2). We give a new proof to this result and propose a kind of converse to get the following result where $GP_nK$ denotes the  set of forms similar to $n$-fold Pfister form  : \vspace *{0.2cm}

\noindent \textbf{Theorem  \ref{decadix}.} \itshape 
Let $ \f $ and $ \psi $ be two anisotropic $ K $-quadratic forms
of dimensions $ \geq 2 $. $ \f $ is
isotropic over $ K (\psi) $, if and only if for 
every field extension $F/K$ and for every form $ \pi \in GP_{n} F $, $ \pi \otimes \f $
is isotropic over $ F (\pi \otimes \psi) $. 
\normalfont \vspace *{0.2cm}

Then, we apply theorem \ref{decadix} to stably birational equivalence to prove  :   \vspace*{0.2cm}

\noindent \textbf{Theorem  \ref{ikakou}.} \itshape Let $ \f $ and $ \psi $ be two anisotropic $ K $-quadratic forms
of dimensions $ \geq 2 $, then the following assertions are
equivalent : 
\begin{enumerate}
\item $ \f \stb \psi $ over $ K $ ; 

\item for every field extension $F/K$,
$ D_{F} (\psi) D_{F} (\psi) = D_{F} (\f) D_{F} (\f) $ ;

  \item for every
$ n \geq1,$  $$\l \l X_{1}, \ldots, X_{n} \r \r \otimes \f \stb \l \l
X_{1}, \ldots, X_{n} \r \r \otimes \psi $$
over $ K (X_{1}, \ldots, X_{n})$ ; 

\item for every field extension $F/K$, for every $ n \geq 1 $ and for
every $ \tau \in GP_{n} (F) $, $$ \f \otimes \tau \stb \psi \otimes \tau.$$
\end{enumerate}
\normalfont \vspace *{0.2cm}

We then apply theorem  \ref{iocentrik} and a work of Rost (\cite{spinor}) to obtain a 
stable birational equivalence criterion for quadratic forms involving an isomorphism of
modified Clifford groups of the considered forms. More precisely, using the notation of \cite{spinor} and the groups $S{\Gamma_{F} (\f)}^{\diamondsuit}$ and $\overline{S \Gamma_{F} (\f)}^{\diamondsuit} $ defined page \pageref{clifgroups}, we get the following result  :  \vspace*{0.2cm}

\noindent \textbf{Theorem  \ref{letheoremsuivantrost}.} \itshape Let $ \f $, $ \psi $ be two anisotropic $ K $-quadratic forms. The following assertions are equivalent : 
\begin{enumerate}
\item $ \f \stb \psi $ ; 

\item for every field extension $F/K$, there exists a group isomorphism
$$ \xi_{F} :  S{\Gamma_{F} (\f)}^{\diamondsuit} \rightarrow
S{\Gamma_{F} (\psi)}^{\diamondsuit} $$ such that for every $ x \in
S{\Gamma_{F} (\f)}^{\diamondsuit}, $ $ \sn^{\diamondsuit} (x) =
\sn^{\diamondsuit} (\xi_{F} (x)) \ ;  \ $

\item for every field extension $F/K$, there exists a group isomorphism
$$ \overline{\xi} _{F} : 
 \overline{S \Gamma_{F} (\f)}^{\diamondsuit} \rightarrow
 \overline{S \Gamma_{F} (\psi)}^{\diamondsuit} $$
 such that for every $ \overline{x} \in
\overline{S \Gamma_{F} (\f)}^{\diamondsuit}, $ $ \overline{\sn
}^{\diamondsuit} (\overline{x}) = \overline{\sn}^{\diamondsuit
} (\overline{\xi} _{F} (\overline{x})). $
\end{enumerate}
\normalfont \vspace *{0.2cm}

Note that, if $\dim(\f)=\dim(\psi)\leq 7$, all those stably birational equivalence conditions are also equivalent to $K(\f)\simeq_K K(\psi)$, as shown in \cite{dim5}, \cite{dim6}, \cite{lag2}, \cite{izhbokar1}, \cite{magicknebusch2}, \cite{magicdetlev}, \cite{DetlevWad}, \cite{lag}, \cite{nr}, \cite{kar} and \cite{moaossijesuiscite}.

\section{Basic definitions and  useful results}

Throughout this paper, we will only consider finite-dimensional nondegenerate
quadratic forms, and we call them forms for short. We
always assume fields to have characteristic different from 2. 
Most of the notations we use are borrowed from the books by Lam (\cite{lam} and \cite{lam2}) and Scharlau (\cite{scharlau}). 
The reader can refer to those books for general background on quadratic forms.

 Let
$\f$ and $\psi$ be forms over $K$. If they are
isometric, we write $\f\simeq\psi$. They are said to be similar
if there exists $a\in K$ such that $\f\simeq a\psi$. We will denote
the orthogonal sum of $\f$ and $\psi$ by $\f\perp\psi$ and their tensor
product by $\f\otimes\psi$. 
We say that $\psi$ is a subform of $\f$  if there exists a form $\sigma$ such that
$\f\simeq\psi\perp\sigma$. In this case,  we write
$\psi\subset\f$.\\
Let $\H$ denote the hyperbolic plane
$\l-1,1\r$.
A form $\f$ is said to be isotropic if
it represents 0 non trivially or, equivalently, if $\H\subset\f$. It
is called hyperbolic if $\f\simeq\H\perp...\perp\H$.
If two forms $\f$ and $\psi$ are Witt-equivalent, we write
$\f\sim\psi$, in particular, $\f\sim0$ means that $\f$ is
hyperbolic.\\
If we consider $\f$ as a form over an extension $F/K$, we write $\f_F$.
$D_K (\f)$ is the set of elements of $K^{\ast}$ represented by
$\f$ and $G_K(\f)$ is the set of similarity factors of $\f$, that is to say, the group of elements $a\in K^\ast$ such that $a\f\simeq \f$.
The multiplicative subgroup
generated by $ D_{K} (\f) $ in
$ K^\ast$ is denoted by $ T_{K} (\f) $, note that this group contains
$ G_{K} (\f) $.\\
Let $\f$ be a  form of dimension $\geq3$, then  the
function field of $\f$ is denoted by $K(\f)$.  In particular, if
$\f\simeq\l a_{1},\ldots,a_{n}\r$, then
$$\displaystyle K(\f)=K(X_{2},\ldots,X_{n})\left(\sqrt{-\frac{1}{a_{1}}(a_{2}X_{2}^{2}
+\ldots+a_{n-1}X_{n-1}^{2}+a_{n}X_{n}^{2})}\right).$$ 
Note that all the results in this paper regarding $K(\f)$ also hold if we consider the projective function filed of $\f$.
Note that for every $ a $ in $ K^{\ast} $, we have
$ K (\f) = K (a \f) $. We will sometimes use that $ K (\f) / K $ is a purely transcendental extension if and
only if $ \f $ is isotropic over $ K $.\\
An $n$-fold Pfister form over $K$ is a form of type
$\l1,-a_1\r\otimes\ldots\otimes\l1,-a_n\r$ where $a_i\in K^{\ast}$, it
will be denoted by $\l\l a_1,\ldots,a_n\r\r$. The set of forms
isometric (resp. similar) to $n$-fold Pfister form will be denoted by $P_n K$ (resp. $GP_n K$).\\
In this paper, $(X_1,\ldots,X_n)$ and $(Y_1,\ldots,Y_n)$ always denote algebraically independent variables.

\section{Analogues of the third representation theorem}

The well known Pfister's third representation theorem is the following (\cite{scharlau}, chapter
4, Theorem 3.7) : 

\begin{theo}[Pfister, 1965]
\label{3rep}
 Let $ \f \simeq \l a_{1}, \ldots, a_{n} \r $ and $ \psi \simeq \l
b_{1}, \ldots, b_{m} \r $ two anisotropic forms over $K$,
then the following assertions are equivalent : 
\begin{enumerate}
\item $ \psi \subset \f $ ;  \item $ D_{F} (\psi) \subset D_{F} (\f) $ for
every field extension $F/K$ ;  \item $ \f $ represents
$ \psi (t_{1}, \ldots, t_{m}) = \overset{n}{\underset{i = 1}{\sum}}b_{i} t_{i}^{2} $ over
$ K (t_{1}, \ldots, t_{m}) $.
\end{enumerate}
In particular, $ \f \simeq \psi $ if and only if for every field extension $F/K$, $ D_{F} (\psi) = D_{F} (\f). $
\end{theo}

This result gives sufficient and necessary conditions equivalent for $ \f $ to be a subform of $ \psi $, and Knebusch's norm theorem, formulated in \cite{norme}, gives an
equivalence between hyperbolicity and a condition of the same type.

\begin{theo}[Knebusch, 1973] \label{gfdthfdsrt}
\index{norm theorem} Let $ \f $ and $ \psi $ be anisotropic forms
such that $ \psi $ represents 1 and $\dim\psi=n\geq2$. The
following assertions are equivalent :  \begin{enumerate}

\item $ \psi (X_{1}, \ldots, X_{n}) \in G_{K (X_{1}, \ldots, X_{n})} (\f) $
 ; 

\item $ \f_{K (\psi)} \sim0 $. \end{enumerate}
\end{theo}

\noindent \textbf{Proof  :  } See \cite{norme}, Theorem 4.2. $ \Box $ \\

We will say that a polynomial is
unitary over the ring $ K [X_{1}, \ldots, X_{m}] $ if its leading coefficient with respect to the lexicographical ordering of monomials 
 is 1. Then, we denote by $ K (g) $ the field of
fractions of $ K [X_{1}, \ldots, X_{m}] / (g) $ for an irreducible polynomial
 $ g $ of $ K [X_{1}, \ldots, X_{m}] $.
Then the well known formulation of Knebusch's norm theorem above is actually a corollary of the more general following result  :  

\begin{theo} [Knebusch, 1973]
Let $ \f $ be an anisotropic quadratic form such that $ 1 \in
D_{K} (\f) $. Let $ a \in K^{\ast} $ and let $ f_{1},\dots f_{r} \in K [X_{1}, \ldots, X_{m}] $ be distinct unitary irreducible polynomials, and let $ f = af_{1} \ldots f_{r} $.

\noindent The following assertions are equivalent : 
\begin{enumerate}
\item $ f \in G_{K (X_{1}, \ldots, X_{m})} (\f) $ ; 

\item $ a \in G_{K} (\f) $ and $ f_{i} \in G_{K (X_{1}, \ldots, X_{m})} (\f) $
for every $ i = 1, \ldots, r $ ; 

\item $ a \in G_{K} (\f) $ and $ \f_{K (f_{i})} $ is hyperbolic for
all $ i = 1, \ldots, r $. \end{enumerate}
\end{theo}

\noindent \textbf{Proof  :  } See \cite{norme}, Theorem 4.2. $ \Box $ \\

Note that if $ \f $ represents 1, the
condition $ \psi (X_{1}, \ldots, X_{n}) \in
G_{K (X_{1}, \ldots, X_{n})} (\f) $ in theorem \ref{gfdthfdsrt} implies
$ \psi (X_{1}, \ldots, X_{n}) \in D_{K (X_{1}, \ldots, X_{n})} (\f) $, just like 
$ \f_{K (\psi)} \sim0 $ implies $ \psi \subset \f $ (using the subform theorem from \cite{wad}, theorem 2).\\
We can then ask if we can exhibit a similar
kind of condition for the isotropy of $ \f $ over $ K (\psi) $.
Note first that, as early as 1956, Witt, in \cite{witt} (see
\cite{toutwitt}) had presented an isotropy criterion
 similar to the third representation theorem (a theorem which was not
yet known). In Witt's formulation, the conditions only concerned irreducible polynomials in one
variable represented by a quadratic form, but Knebusch noticed in \cite{norme} that the argument was valid
for irreducible polynomials in $ m $ variables.
We recall that
$ T_{K (X_{1}, \ldots, X_{m})} (\f) $ is the multiplicative subgroup
generated by $ D_{K (X_{1}, \ldots, X_{m})} (\f) $ in
$ K (X_{1}, \ldots, X_{m}) $. Bayer-Fluckiger 
developed Knebusch's proof to formulate in
\cite{eva} the following theorem :

\begin{theo} [Witt, 1956  ;  Knebusch, 1973  ;  Bayer-Fluckiger, 1995]
\label{faich} \hfill

\noindent Let $ \f $ be an anisotropic quadratic form which
represents 1, let $ a \in K^{\ast} $, let $ f_{1},\dots f_{r} \in K [X_{1}, \ldots, X_{m}] $ be distinct unitary irreducible polynomials, and let $ f = af_{1} \ldots f_{r} $. 

\noindent The following assertions are equivalent : 
\begin{enumerate}
\item $ f \in T_{K (X_{1}, \ldots, X_{m})} (\f) $ ; 

\item $ a \in T_{K} (\f) $ and $ f_{i} \in T_{K (X_{1}, \ldots, X_{m})} (\f) $
for every $ i = 1, \ldots, r $ ; 

\item $ a \in T_{K} (\f) $ and $ \f_{K (f_{i})} $ is isotropic for every
$ i = 1, \ldots, r $. \end{enumerate}
\end{theo} 

\noindent \textbf{Proof  :  } See \cite{eva}, Theorem 1. $ \Box $ \\

In particular, if $ r = 1 $ and  $ f $ is a
quadratic form, we get the following characterizations
  : 

\begin{theo} \label{letheormequiresumelesautres} Let $ \f $ and $ \psi $ be two anisotropic  $ K $-quadratic forms with $\dim\psi\geq2$. Suppose
that $ \psi $ and $ \f $ represent 1, then with the notation
$ X = X_{1}, \ldots, X _{\dim (\psi)} $,
\begin{enumerate}
\item \textbf{\textup{(Knebusch, 1973).}}\\ $ \f $ is hyperbolic
over $ K (\psi) $ if and only if $ \psi (X) \in G_{K (X)} (\f) $ ; 

\item \textbf{\textup{(Pfister, 1965).}}\\ $ \psi $ is a subform of $ \f $ if
and only if $ \psi (X) \in D_{K (X)} (\f) $ ; 

\item \textbf{\textup{(Witt, 1956 ;  Knebusch, 1973 ; Bayer-Fluckiger, 1995).}}\\ $ \f $ is isotropic over $ K (\psi) $ if and
only if $ \psi (X) \in T_{K (X)} (\f) $. \end{enumerate}
\end{theo}

\noindent \textbf{Proof  :  } See \cite{eva}, Theorem 9. $ \Box $ \\

\subsection{An alternative solution}

We are going to reformulate theorem  \ref{faich} in the case of
quadratic forms and we will consider sets that are in general smaller than $ T_{K} (\f) $, namely $ D_{K} (\f) D_{K} (\f) $ and the multiplicative group
$ N_{K} (\f) $ generated
 by $ D_{K} (\f) D_{K} (\f) $ in
$ T_{K} (\f) $.\\ 
Note that $ N_{K} (\f) $ is the group of norms of a
quadric. Indeed, let $ X $ be a variety over $ K $, then the group
of norms $ N_{K} (X) $ is generated by all norms $ z\in
N_{L / K} (L ^{\ast}) $ where $ L $ runs over every finite extension of
$ K $ for which $ X $ admits an $ L $-rational point. If $ X = Q $ is the
quadric associated with a nonsingular quadratic form $ q $ of
dimension $\geq 3 $, then Colliot-Thélène and Skorobogatov showed in lemma  2.2 of \cite{colliotskoro}
that $ N_{K} (Q) = N_{K} (q) $.

\begin{rem}Note that we have $[T_K(\f):N_K(\f)]=1\ or\ 2$, and $T_K(\f)=N_K(\f)$ if $1\in D_K(\f)$.\\
For example, if $\f=\l a\r$ with $a\in K^\ast\setminus K^{\ast2}$, then $T_K(\f)=K^{\ast2}\cup aK^{\ast2}$ whereas
$N_K(\f)=K^{\ast2}$, hence $[T_K(\f):N_K(\f)]=2$. However, if, say, $K=K_0((X))((Y))$ is the iterated Laurent series field in two variables over a field $K_0$ of characteristic not 2, and if $\f=\l X,Y,XY\r$, then $1\not\in D_K(\f)$.
Nonetheless, $T_K(\f)=N_K(\f)=K^{\ast2}\cup XK^{\ast2}\cup YK^{\ast2}\cup XYK^{\ast2}=D_K(\l\l -X,-Y\r\r)$.
\end{rem}

\begin{theo}
\label{iocentrik} Let $ \f $ and $ \psi $ be two anisotropic $ K $-quadratic forms with $\dim\psi\geq2$.\\ If we write
$ X = X_{1}, \ldots, X _{\dim (\psi)} $, then the following assertions
are equivalent : 
\begin{enumerate}
\item $ \f_{K (\psi)} $ is isotropic ; 

\item $ D_{F} (\psi) D_{F} (\psi) \subset D_{F} (\f) D_{F} (\f) $ for
every field extension $F/K$ ; 

\item $ a \psi (X) \in
D_{K (X)} (\f) D_{K (X)} (\f) $ for
every $ a \in D_{K} (\psi) $ ; 

\item $ N_{F} (\psi) \subset N_{F} (\f) $ for
every field extension $F/K$ ;

\item $ a \psi (X) \in N_{K (X)} (\f) $
for every $ a \in D_{K} (\psi) $ ;

\item $ T_{F} (a \psi) \subset T_{F} (\f) $ for 
every field extension $F/K$  and every $ a \in D_{K} (\psi) $ ;

\item $ a \psi (X) \in T_{K (X)} (\f) $
for every $ a \in D_{K} (\psi) $.
\end{enumerate}
\end{theo} 

\noindent The proof of this result will use the four following lemmas.

\begin{lem}\label{pitilemme}For every $ a\in D_{K} (\f) $, $ N_{K} (\f) = T_{K} (a\f)$.
\end{lem}

\noindent\textbf{Proof  : }  Let $V$ be the underlying vector space of $\f$. If $ a =\f (y) $ for some $y\in V$, and if $ x\in T_{K} (a\f) $, then $ x =\overset{r}{\underset{i=1}{\prod}}a\f (x_{i}) $ for some $x_1,\ldots x_r \in V$ and we have $ x =\overset{r}{\underset{i=1}{\prod}}\f (y)\f (x_{i})\in N_{K} (\f) $.\\ Conversely, let $ x\in N_{K} (\f) $,
then $ x =\overset{r}{\underset{i=1}{\prod}}\f (y_{i})\f (x_{i}) =\overset{r}{\underset{i=1}{\prod}}\f (y_{i})\overset{r}{\underset{i=1}{\prod}}\f (x_{i}) $ for some $x_i,y_i \in V$, $1\leq i \leq r$,
and thus, we will have $ x =\overset{r}{\underset{i=1}{\prod}} a\f (\frac{y_{i}}{a})\overset{r}{\underset{i=1}{\prod}}
a\f (x_{i}) $. Therefore $ x\in T_{K} (a\f) $.$\Box $

\begin{lem}\label{lemdetransitivite}Let $\f$, $\psi$, $\xi$ be quadratic forms over $K$ with $\psi$ isotropic over $K(\xi)$. If $\f$ is isotropic over $K(\psi)$, then $\f$ is isotropic over $K(\xi)$
\end{lem}

\noindent\textbf{Proof  : } If $\f$ is isotropic over $K(\psi)$, then it is so over $K(\psi)(\xi)=K(\xi)(\psi)$. But $\psi$ is isotropic over $K(\xi)$, thus $K(\xi)(\psi)$ is purely transcendental over $K(\xi)$ and hence $\f$ is already isotropic over $K(\xi)$.$\Box$

\begin{lem}\label{pouet} Let $\f $ and $\psi $ be two
$ K $-forms of dimensions $\geq 2 $. If $\f $ is isotropic over $ K (\psi) $, then for
every field extension $F/K$,
$ D_{F} (\psi) D_{F} (\psi)\subset D_{F} (\f) D_{F} (\f). $\end{lem}

\noindent\textbf{Proof  : }  If $\f $ is isotropic over $ K $,
it will be so over any extension $ F $
 over $ K $. Therefore $ D_{F} (\f) D_{F} (\f) = F^{\ast} $ for every field extension $F/K$ and the
property is verified. If $\psi $ is isotropic, then $ K (\psi) $ is a purely transcendental extension
 of $ K $. Since $\f $ is a $ K $-form isotropic 
over a purely transcendental extension of $ K $, it is isotropic already over
$ K $. According to the previous case, the property is still verified.\\
Suppose that $\f $ is anisotropic over $ K $ and isotropic over
$ K (\psi) $, then for every  field extension $F/K$, $\f $ will also be
isotropic over $ F (\psi) $. We may also assume that $\f $ is anisotropic
over $ F $, otherwise the property is obvious.
Let then $ c\in D_{F} (\psi) $, we can write $\psi\simeq\l
c\r\perp c\tau $ and therefore the sets $ cD_{F} (\psi) $,
$ D_{F} (c\psi) $ and $ D_{F} (\l1\r\perp\tau) $ are equal.
 Consider
now an element $ a $ of $ D_{F} (c\psi) $, we can write
$ a = x^{2} + y $ with $x\in F$ and $ y\in D_{F} (\tau) $ or $ y = 0 $.\\
 If $ y\not = 0 $, then
$ c\psi\simeq\l1, y\r\perp\eta $. But then $\psi_{F (\l1, y\r)} $ is
isotropic and  $\f_{F (\l1, y\r)} $ is isotropic by lemma \ref{lemdetransitivite}, which
gives,
using lemma 4.2 from \cite{detlevsurvey}, $\f\simeq b\l1, y\r\perp\xi $ for some $ b\in F $. As a consequence, 
we will have $ x^{2} + y = (bx^{2} + by)\frac{b}{b^{2}}\in
D_{F} (\f) D_{F} (\f). $\\ 
If $ y = 0 $, for every $ b\in D_{F} (\f) $,
$ a = x^{2} = bx^{2}\frac{b}{b^{2}}\in D_{F} (\f) D_{F} (\f). $ For every
$ c\in D_{F} (\psi) $, we thus get that $ cD_{F} (\psi)\subset
D_{F} (\f) D_{F} (\f). $\\
Then, if $ d\in D_{F} (\psi) D_{F} (\psi) $, we can
write $ d = ef $ where $ e, f\in D_{F} (\psi) $. But then it is clear
that $ d\in eD_{F} (\psi) $ and therefore, from the above, $ d\in
D_{F} (\f) D_{F} (\f) $. $\Box $

\begin{lem}\label{lemmepourfaireplusselfcontained}Let $f\in K[X_1,\ldots,X_m]=K[X]$ be an irreducible polynomial and $\f$ a quadratic form over $K$.
If there exists $a\in K^\ast$ such that $af\in T_{K(X)}(\f)$, then $\f_{K(f)}$ is isotropic.
\end{lem}

\noindent \textbf{Proof  : }Say $\dim\f=n$. Since replacing $f$ by $af$ does not change $K(f)$, we may assume that $a=1$.
Write $f=\overset{r}{\underset{i=1}{\prod}}\f\left(\frac{P_{i1}}{g_i},\ldots,\frac{P_{in}}{g_i}\right)$, with $P_{i,j}, g_i\in K[X]$
and $g_i\not=0$. Let $g=\overset{r}{\underset{i=1}{\prod}}g_i$. Then $fg^2=\overset{r}{\underset{i=1}{\prod}}\f(P_{i1},\ldots,P_{in})$.
If there exists $i$ such that $\forall 1\leq j\leq n,\ f|P_{ij}$, then $f^2|fg^2$.
Then $f|g$ since $f$ is irreducible and we may replace $P_{ij}$ by $\frac{P_{ij}}{f}\in K[X]$ and $g$ by $\frac{g}{f}\in K[X]$.
Repeating this if necessary, we may eventually assume that to each $i\in\{1,\ldots,r\}$ there exists $j\in\{1,\ldots,n\}$ such that $f\not|P_{ij}$.\\
Since $f|\overset{r}{\underset{i=1}{\prod}}\f(P_{i1},\ldots,P_{in})$ and $f$ irreducible, there exists $i\in\{1,\ldots,r\}$ such that
$f|\f(P_{i1},\ldots,P_{in})$.\\
Then there exists $h\in K[X]$ such that $\f(P_{i1},\ldots,P_{in})=fh$.\\
Since not all $P_{ij}$ divisible by $f$, this shows that $\f$ is isotropic over the integral domain $K[X]/(f)$ and hence also over its quotient field $K(f)$.$\Box$\\

\noindent \textbf{Proof of theorem  \ref{iocentrik} : } First, $ (1) \Rightarrow
(2) $ is given by lemma   \ref{pouet}.\\
We  will now use that, by definition,
$ D_{F} (\f) D_{F} (\f) \subset N_{F} (\f) \subset T_{F} (\f) $ 
for every field extension $F/K$ . This shows in particular that
 $ (3) \Rightarrow(5) \Rightarrow(7) $.
Moreover, since $ N_{F} (\f) $ is the multiplicative group generated
by $ D_{F} (\f) D_{F} (\f) $, we will always have $ (2) \Rightarrow(4) $. Finally,
using lemma \ref{pitilemme}, we have $ T_{F} (a \psi) = N_{F} (\psi) $ for all $a\in D_K(\psi)$. Then, as $ N_{F} (\f) \subset T_{F} (\f) $
for
every field extension $F/K$ , we will have $ (4) \Rightarrow(6) $.
 As $ (3) $ and $ (7) $ are respectively special cases of $ (2) $ and
 $ (6) $, we will have $ (2) \Rightarrow(3) $ and $ (6) \Rightarrow(7) $.
Thus, we get
$ (1) \Rightarrow(2) \Rightarrow(4) \Rightarrow(6) \Rightarrow(7) $
 and $ (1) \Rightarrow(2) \Rightarrow(3) \Rightarrow(5) \Rightarrow(7) $.\\
Since $\psi(X)$ is irreducible, lemma \ref{lemmepourfaireplusselfcontained} shows that $(7) \Rightarrow(1)$.$\Box$
%Let us show that $ (7) \Rightarrow (1) $\\
% Since replacing $ \psi $ by
%$ a \psi $ does not change anything regarding isotropy,
%we may assume that $ a = $ 1. While $ \psi \simeq \l a_{1}, \ldots \r $
%and $ p = \psi (X_{1}, \ldots, X_{n}) $, and consider the field
%$ K (X_{1}, \ldots, X_{n}) (b) $ with
%$ b = \sqrt{-\frac{1}{a_{1}} (\psi (0, X_{2}, \ldots, X_{n}))}. $ The field
%$ K (X_{1}, \ldots, X_{n}) (b) $ is the same as
%$ K (X_{2}, \ldots, X_{n}) (b) (X_{1}) $. But over this field, we have 
%$ p = a_{1} X_{1}^{2} -a_{1} b^{2} $ and hence
%$ p = a_{1} X_{1}^{2} -a_{1} b^{2} = a_{1} (X_{1} -b) (X_{1} + b). $
%Then, for $ Y = X_{1} -b $ it is
%clear that $ K (X_{2}, \ldots, X_{n}) (b) (Y)= K (X_{2}, \ldots, X_{n}) (b) (X_{1}) $. But on the other hand,
% $ K (X_{2}, \ldots, X_{n}) (b) (Y) \simeq K (\psi)^{\#} (Y) $ and over
%this field, $ p = a_{1} Y (Y + 2b) $ is an element of valuation 1 for the $ Y $-adic
%valuation  over $ K (\psi)^{\#} (Y) $
% belonging
%to the group $ T_{K (X_{1}, \ldots, X_{n})} (\f) $ and thus to the group
%$ T_{K (\psi)^{\#} (Y)} (\f) $.\\ 
%As a consequence, $\f$ represents an element of valuation 1 over $ K (\psi)^{\#}(Y)$.
%Lemma   \ref{lanouvellepropquiresoiuttout} now implies that $ \f $ is
%isotropic over $ K (\psi)^{\#} $ and therefore over $ K (\psi)$.$ \Box $

\subsection{More representation results}

In her proof of the theorem  \ref{faich}, 
Bayer-Fluckiger uses the following result, a variant of which we will prove below : 

\begin{lem} Let $ \f $ be a quadratic form of dimension $ n $ and let $ f \in K [X_1,\ldots,X_m] $ be an irreducible polynomial such that its leading coefficient with respect to the lexicographical ordering is in $N_K(\f)$. If
 $ f $ divides a polynomial of the form
$ \f (p_{1}, \ldots, p_{n}) $ with $ p_{i} \in K [X_1,\ldots,X_m] $ setwise coprime, then $ f\in N_{K (X_1,\ldots,X_m)}(\f) $.
\end{lem}

\noindent\textbf{Proof : } \cite{eva} \S 3, lemme.$\Box$\\

\noindent It seems interesting to formulate the following variant because it
leads to consider a simpler formulation of the third point of theorem
 \ref{letheormequiresumelesautres}. The set $\left\lbrace \overset{r}{\underset{i=1}{\prod}}a_i|a_i\in D_K(\f)\right\rbrace$ of products or $r$ elements of
$ D_{K} (\f) $ will be denoted by $D_K(\f)^r$.

\begin{pro}
\label{trh} Let $ f \in K [X] $ be an irreducible unitary polynomial in
one variable $ X $. Let $ \f $ be a quadratic form of dimension $ n $. If
 $ f $ divides a polynomial of the form
$ \f (p_{1}, \ldots, p_{n}) $ with $ p_{i} \in K [X] $ not all divisible
by $ f $, then $ f\in D_{K (X)} (\f)^m $  with $ m \leq \deg_{X} (f) $.
\end{pro}

\noindent \textbf{Proof  :  } If $ \dim (\f) = $ 1 then $ \f = \l
a \r $ and $ ap^ 2 = fg $, so $ f $ divides $ p $, which contradicts
the hypothesis that the $ p_i $ are not all divisible by $ f $. The
case of dimension 1 is therefore excluded. If $ \f $ is isotropic, then we
have $ f \in D_{K (X)} (\f) = K (X)^{\ast} $. We may therefore assume that $ \f $
is anisotropic of dimension $ \geq2 $. Let us consider the degree
of $ f $.\\
If $ \deg_{X} (f) = 1 $, then there exists $ a \in K $ such that $ f = X-a $.
Suppose then that $ (X-a) g (X) = \f (p_{1}, \ldots, p_n) $. In this
case, $ 0 = \f (p_{1} (a), \ldots, p_{n} (a)) $. As $ \f $ is
anisotropic, we will have $ p_i (a) = 0 $ for every $ i $. But then every
$ p_i $ is divisible by $ f = X-a $, which contradicts our
hypotheses.\\
If $ \deg_{X} (f) = 2 $, after a linear change of variable, we
may assume that there exists $ a \in K^* $ such that $ f =
X^ 2-a $,  an irreducible polynomial. We can then write $ (X^ 2-a) g (X) =
\f (p_{1}, \ldots, p_n) $. Let us consider $ K (\sqrt{a}) $ and
evaluate in $ X = \sqrt{a} $. If $ p_{i} (\sqrt{a}) = 0 $ for every $ i $,
then there exists $ q_{i} \in K [X] $ such that $ p_{i} = (X^{2} -a) q_{i} $
for every $ i $, which contradicts our hypothesis. So we may 
assume that there exists $ k $ such that $ p_{k} (\sqrt{a}) \not = 0 $. As
$ \f (p_{1} (\sqrt{a}), \ldots, p_n (\sqrt{a})) = 0 $ with the
$ p_{i} (\sqrt{a}) $ not all zero, it follows that $ \f $ becomes
isotropic over $ K (\sqrt{a}) $. From \cite{scharlau}, chapter 2, lemma  5.1, we will then have
$ \f \simeq c \l 1, -a \r \perp \f '$. So $ \f $ represents
$ \frac{1}{c} $ and $ c (X^ 2-a) $ over $ K (X) $, so $ f = X^ 2-a \in
D_{K (X)} (\f) D_{K (X)} (\f). $ \\
Let us then proceed by induction on the degree. We can divide the
$ p_{i} $ by $ f $ to get $ p_{i} = \xi_{i} f + \overline{p_{i}} $ where
$ \deg_{X} (\overline{p_{i}}) <\deg_{X} (f) $. Recall that, if
$ p = (p_{1}, \ldots, p_{n}) $,  then we will have by hypothesis $ \f (p) = fg $ for a
certain polynomial $ g $. In this case, if
$ \xi = (\xi_{1}, \ldots, \xi_{n}) $ and
$ \overline{p} = (\overline{p_{1}}, \ldots, \overline{p_{n}}) $, then we 
can write : 
\begin{eqnarray*}
\f (\overline{p}) & = & \f (p) +
f^{2} \f (\xi) -2fb _{\f} (p, \xi) \\
 & = & f (\underbrace{g +
f \f (\xi) -2b _{\f} (p, \xi)} _{ah}) \\
& = & afh.
\end{eqnarray*}
with $ a \in K^{\ast} $ and $ h \in K [X] $ unitary. Consider
$ h = \prod h_{j} $ where the $ h_{j} \in K [X] $ are irreducible
unitary, then
$ \deg_{X} (f) + \sum \deg_{X} (h_{j}) \leq
2 \max_{i} (\deg_{X} (\overline{p_{i}})) <2 \deg_{X} (f). $ As a 
consequence,
$ \deg_{X} (h_{j}) \leq
\underbrace{\sum \deg_{X} (h_{j})} _{\deg_{X} (h)} <\deg_{X} (f). $ 

\noindent Since not all $ \overline{p_{i}}=0 $ and since $\f$ is anisotropic, a simple degree argument yields that $\deg_{X} (\f (\overline{p}))$
is even, hence $ \deg_{X} (h) \leq
\deg_{X} (f) -2 $.
If an $ h_{k} $ divides every $ \overline{p_{i}} $ then $ h_{k}^ 2 $ divides $\f(\overline{p})=afh$, and thus $ h_{k}^ 2 $
divides $h$ since both $f$ and $h_k$ are irreducible and $ \deg_{X} (h_{k}) <\deg_{X} (f) $.
After cancelling such an $ h_{k}^ 2 $ from $h$ and $h_k$ from each $\overline{p_i}$, we may in fact assume that we have 
$afh=af\underset{j}{\prod}h_j=\f(p')$ with $ p '= (p_{1}', \ldots, p_{n} ') \in K[X]^n$, and the $h_j$ irreducible
unitary and not dividing every $ p_{i} '$.\\
Let us first show that $ a \in D_{K} (\f) $. As $ f $ and the $ h_{j} $ are
unitary, the highest degree coefficient of $ af \prod h_{j} $
is exactly $ a $. Let's denote by $ a_{i} '$ the highest degree coefficient
 of each $ p_{i} '$ and let $ d $ be the maximum degree of
$ p_{i} '$, we may then write
$ \widetilde{a_{i}} = \left \lbrace \begin{array}{l} a_{i} '\
\mathrm{if} \ \deg_{X} (p_{i} ') = d \\
0 \ \mathrm{otherwise} \end{array} \right. $ and
$ \widetilde{a} = (\widetilde{a_{1}}, \ldots, \widetilde{a_{n}}) $.
Since $ \f $ is anisotropic, the coefficient of highest degree of
$ \f (p ') $ is $ \f (\widetilde{a}) $. Thus, $ a = \f (\widetilde{a}) $ and
$ a \in D_{K} (\f) $.\\
If $h=1$, then $ f = \frac{1}{a} \f (p')\in D_{K} (\f) D_{K (X)} (\f)\subset D_{K (X)}(\f)^2$ with $2\leq \deg_{X} (f)$.\\
If $h\not=1$, then as $ \deg_{X} (h_{j}) <\deg_{X} (f) $ we can apply the induction hypothesis which implies that
$h_{j} \in D_{K (X)} (\f)^{m_j}$ with $m_j\leq\deg_{X} (h_{j})$, hence with $h=\prod h_{j}$ we get 
$ f = \frac{1}{a h} \f (p ') \in D_K(\f)D_{K (X)} (\f)^{\sum m_j}D_{K (X)} (\f)$ and thus $f\in D_{K (X)} (\f)^{\sum m_j + 2}$
with $\sum m_j + 2\leq\deg_{X} (h)\leq\deg_{X} (f) $.$ \Box $ \\

We can also use the proposition  \ref{trh} to show
a more general property than the equivalence $ 1 \Leftrightarrow3 $ of the theorem
 \ref{iocentrik}, and more precise than the theorem  \ref{faich} in
the case of polynomials in one variable.

\begin{pro} \label{etjevousdispouet} Let $ \f $ be an anisotropic quadratic form
representing 1. Let $ f \in K [X] $ be an irreducible unitary polynomial in one variable. Then the following assertions are equivalent : 

\begin{enumerate}
\item $ f \in D_{K (X)} (\f)^m$ with $m\leq \deg (f)$ ; 

\item $ \f_{K (f)} $ is isotropic. \end{enumerate}
\end{pro}

\noindent The proof will need the following result.

\begin{lem} \label{encoreunpetitlemmedeplus} If a quadratic form
 $ \f $ of dimension $n$ over $K$ becomes isotropic over $ K (f) $ where $ f \in
K [X_1,\ldots, X_m]=K[X] $ is an irreducible polynomial, then $ f $ divides
a polynomial of the form $ \f (p_{1}, \ldots, p_{n}) $ with the
$ p_{i} \in K [X_1\ldots X_m] $ not all divisible by $ f $.
\end{lem}

\noindent \textbf{Proof  : }Since $\f$ is isotropic over 
$K(f)=\mathrm{Quot}(K[X]/(f))$, after clearing denominators we can find $\overline{p_1},\ldots,\overline{p_n}\in K[X]/(f)$
not all $\overline{p_i}=0$ such that $\f(\overline{p_1},\ldots,\overline{p_n})=0$ in $K[X]/(f)$, or, equivalently,
we can find $p_1,\ldots,p_n,h\in K[X]$, not every $p_i$ divisible by $f$, such that $\f(p_1,\ldots,p_n)=fh$.$\Box $\\

\noindent \textbf{Proof of the proposition
 \ref{etjevousdispouet} : } $ (1) \Rightarrow(2) $ As $ f \in D_{K(X)(\f)^m}$ with $m\leq\deg_X(f)$, this implies that $ f \in
T_{K (X)} (\f) $ and it suffices to apply lemma \ref{lemmepourfaireplusselfcontained}.\\
$ (2) \Rightarrow(1) $ As $ \f_{K (f)} $ is isotropic, according to
lemma   \ref{encoreunpetitlemmedeplus}, $ f $ divides a
polynomial of the form $ \f (p_{1}, \ldots, p_{n}) $ with the $ p_{i} \in
K [X] $ not all divisible by $ f $. It is then sufficient to apply
proposition  \ref{trh}. $ \Box $ \\

In Corollary 1 of \cite{eva}, Bayer-Fluckiger establishes that
if $ \f $ is a form that represents 1 and $ f \in
K [X_{1}, \ldots, X_{n}] $ is an irreducible and unitary polynomial, $ \f $
is isotopic over $ K (f) $ if and only if $ f \in
T_{K (X_{1}, \ldots, X_{n})} (\f) $. We will now present
a finer version of this equivalence where a bound appears
for the number of factors in $ T_{K (X_{1}, \ldots, X_{n})} (\f) $. This
result, which also refines the theorem  \ref{faich} and
implies $ (1) \Leftrightarrow(3) $ in the theorem  \ref{iocentrik}, 
was obtained in collaboration with Detlev W. Hoffmann. 

\begin{theo} \label{coroeva} Let $ \f $ be an anisotropic quadratic form of dimension $n$ which
represents 1. We write $X=X_1,\ldots,X_n$. Let $ f \in K [X] $ be an irreducible polynomial
 where all $ X_{i} $ appear. If the degree of the term of highest degree for the lexicographical order
 is $ M $ and if $ f $ is unitary for the lexicographical order, then the following assertions are equivalent
 :  \begin{enumerate} \item  $ f \in D_{K(X)(\f)^m}$ with $m\leq M$ ; 
\item $ \f_{K (f)} $ is isotropic. \end{enumerate}
\end{theo}

\noindent We will need the following lemma to prove this result.

\begin{lem}\label{lelemmebonusdeDetlev}Let $\f$ be an anisotropic quadratic form of dimension $n$ over $K$, and let $a\in D_{K((Y))}(\f)$.
Let $\nu_Y : K((Y))\rightarrow \Z\cup\infty$ be the $Y-$adic valuation. We write $c=c_1,\ldots,c_n\in K((Y))$ and $a=\f(c)$.\\
Then $\nu_Y(a)=2\min\{\nu_Y(c_i)|1\leq i\leq n\}$.\\
In particular, if $a\in  D_{K((Y))}(\f)^m\cap K$, then $a\in D_K (\f)^m$.
\end{lem}

\noindent\textbf{Proof  :  }The fact that $\nu_Y(a)=2\min\{\nu_Y(c_i)|1\leq i\leq n\}$ follows easily from the anisotropy of $\f$ with a standard argument comparing terms of lowest degree in the $X_i\in K((Y))$.\\
If $a\in  D_{K((Y))}(\f)^m\cap K$, let $a_i\in D_{K((Y))}(\f)$ with $a=\overset{m}{\underset{i=1}{\prod}}a_i$. By the above, there are 
$k_i\in\Z$ such that $\nu_Y(a_iY^{2k_i})=0$.
Put $a_i'=a_iY^{2k_i}\in  D_{K((Y))}(\f)$, $k=\overset{m}{\underset{i=1}{\sum}}2k_i$.
Then $\overset{m}{\underset{i=1}{\prod}}a_i'=aY^k$ and $0=\nu_Y\left(\overset{m}{\underset{i=1}{\prod}}a_i'\right)=\nu_Y(a)+k=k$ since 
$\nu_Y(a)=0$. Hence $k=0$ and $a=\overset{m}{\underset{i=1}{\prod}}a_i'$. Also by the above, we can now write $a_i'=\f(c_{i1},\ldots,c_{in})$
with $c_{ij}\in K((Y))$ and $\nu_Y(c_{ij})\geq0$, that is to say $c_{ij}\in K[[Y]]$.
Passing to the residue field $K$ (that is to say putting $Y=0$) we get $a=a(0)=\overset{m}{\underset{i=1}{\prod}}a_i'(0)$ with 
$a_{i}'(0)=\f(c_{i1},\ldots,c_{in}(0))\in D_K(\f)$, hence the claim.$\Box$\\

\noindent \textbf{Proof of theorem \ref{coroeva}  :  } We will write  $ X '= X_{2}, \ldots, X_{n} $.\\
For $ (1) \Rightarrow (2) $ comes from lemma \ref{lemmepourfaireplusselfcontained}.\\
For $ (2) \Rightarrow(1) $, we proceed by induction on the number of
variables,  the case of one variable is treated by proposition
 \ref{etjevousdispouet}.\\ Let $ m = \deg_{X_{1}} (f) $ and 
$ f = X_{1}^{m} f_{m} +
\underbrace{\ldots \ldots \ldots} _{\deg_{X_{1}} <m} $ 
 where $ f_{m} \in K [X '] $.\\
By Gauss' lemma $F = \frac{f}{f_{m}}\in  K (X ') [X_{1}]$ is a unitary irreducible polynomial in $X_1$ over $K(X')$ of degree $m$, and the degree $f_m$ for the lexicographical order is $M-m$. Note also that $K(X')(F)\simeq K(f)$, so by proposition \ref{etjevousdispouet}, $F$ can be written as a product of at most $ m $ elements in $D_{K(X)} (\f) $ and thus of $m$ elements since $ 1 \in D_K (\f) $.
Moreover, by multiplying $ F $ by $ G^{2} $ for a suitable $ G \in
K [X] $, we may assume that $ FG^{2} $ is a product of
$ m $ elements of $ D_{K [X]} (\f) $, that is 
$ FG^{2} = \f (p_{11}, \ldots, p_{1n}) \ldots \f (p_{m1}, \ldots, p_{mn}) $
with $ p_{ij} \in K [X] $. 
Moreover, we may assume that
if $ g $ is an irreducible factor of $ G $, then there does not exist
any $ k $ such that $ g $ divides every $ p_{k1}, \ldots, p_{kn} $, for otherwise
we may divide every $ p_{kj} $  by $ g $ and $ G^{2} $  by $ g^{2} $.\\
If $ \deg (f_{m}) = 0 $, then since $ f $ is unitary, $ f_{m} = 1 $ and the
proof is complete.\\
If $ \deg (f_{m})> 0 $, then %, since $ \f $ is isotropic over $ K (f) $, it is so also over the integral domain $ K [X] / (f) $ and therefore 
using lemma \ref{encoreunpetitlemmedeplus} there
exists $p_{1}, \ldots, p_{n}\in K [X] $, not all $p_i$ divisible by $f$ and $ h \in K [X] $ such that
$ hf = \f (p_{1}, \ldots, p_{n}) $. But, since $ \f $ is anisotropic, by comparing the coefficients of the highest degree terms in $ X_{1} $, we
obtain   polynomials $ h ', p_{i}' \in K [X '] $ such that
$ h'f_{m} = \f (p '_{1}, \ldots, p' _{n}). $ Then we may write
$ f_{m} = r^{2} s $ where $ r, s \in K [X '] $ are unitary and $ s $ has no
square factors.\\
If $ s = 1 $, then $ f_{m} = r^{2} \in D_{K (X)} (\f) $ and therefore
$ f = f_{m} \frac{f}{f_{m}} $ is a product of at most $ m + 1 \leq
m + \deg (f_m) = deg (f) = M $ elements of $ D_{K (X)} (\f) $, as desired.\\
If $ s \not = 1 $, let $ t \in K [X '] $ be an irreducible factor of $ s $.
Then by considering $ f $ as a polynomial in $ X_{1} $ with
coefficients in $ K [X '] $, the polynomial $ t $ does not divide all
coefficients of $ f $ since $ f $ is irreducible.\\
Suppose that $ K $ is infinite, then there exists $ c \in K $ such that
$ t = t (X ') $ does not divide $ f (c, X') \in K [X '] $. In particular, $f(c,X')\not=0$. But as $ G^{2} F $
is a product of $ m $ elements of $ D_{K [X]} (\f) $ and $ f_{m} \in
K [X '] $, then $ G^{2} ff_{m} = G^{2} f_{m}^{2} F $ is always a
product of $ m $ elements of $ D_{K [X]} (\f) $ and we can therefore
substitute $ X_{1} = c $. So $ G (c, X ')^ 2f_m (X') f (c, X ') $ is always
a product of $ m $ elements of $ D_{K [X ']} (\f) \cup \{0 \} $.\\
If $ G (c, X ') = 0 $, then $ X_{1} -c $ divides $ G $ and the anisotropy
of $ \f $ over $ K (X ') $ implies that for at least one $ i $, we
have $ p_{i1} (c, X ') = \ldots = p_{in} (c, X') = 0. $ Therefore $ X_{1} -c $
divides $ p_{i1}, \ldots, p_{in} $ and $ G $, which contradicts the
properties of $ G $. \\
Thus, $ G (c, X ') \not = 0 $ and therefore $ G (c, X') f_m (X ') f (c, X') \neq 0. $\\
So let us write $ s = s't $ with $ s '\in K [X'] $, and 
$ f '(X') = s' (X ') f (c, X') \in K [X '] $, and $ l (X') = G (c, X ') r (X') \in K [X '] $.
Then we have
\begin{eqnarray*} G (c, X ')^{2} f_{m} (X') f (c, X ') & = &
G (c, X ')^{2} r^{2} (X') s (X ') f (c, X') \\& = & G (c, X ')^{2} r^{ 2} (X ')
s' (X ') t (X ') f (c, X') \\& = & l^{2} f't \end{eqnarray*} with $ t $
not dividing $ f '$.\\ Let us set $ p '_{ij} = p_{ij} (c, X') $. We then have
$ l^{2} f't = \f (p '_{11}, \ldots, p' _{1n}) \ldots \f
(p '_{m1}, \ldots, p' _{mn}). $ If there is an $ i $ such that every
$ p '_{i1}, \ldots, p' _{in} $ is divisible by $ t $, then
$ t^{2} $ divides $ l^{2} f't $. Therefore, $ t $ divides $ l $, and we can
divide $ l $ and the $ p '_{ij} $ by $ t $ and we conclude by an induction on the number of irreducible factors of $s$. 
We can therefore assume that
there does not exist any $i$ such that every $ p '_{i1}, \ldots, p' _{in} $ is divisible by $ t $. 
But there exists $ k $ such that $ t $ divides
$ \f (p '_{k1}, \ldots, p' _{kn}) $, so $ \f $ becomes isotropic over
 the integral domain $ K [X '] / (t) $ and
therefore also over its quotients field $ K (t) $. By induction hypothesis, since $ t $ is irreducible
and unitary, $ t $ can be written as a product of at most $ \deg (t) $
factors in $ D_{K (X ')} (\f) $.\\
By repeating this process for every irreducible factor of $ s $,
 it follows from adding the degrees of the factors of $ s $ that $ s $ and
so $ f_m = r^ 2s $ can be written as a product of at most
$ \deg (s) \leq \deg (f_m) = M-m $ elements of $ D_{K (X ')} (\f) $. Then,
since $ F $ is a product of at most $ m $ elements of $ D_{K (X)} (\f) $,
we can write $ f = Ff_{m} $ as a product of at most
$ \deg (F) + \deg (f_m) = M-m + m = M $ elements of $ D_{K (X)} (\f) $.\\
If $ K $ is a finite field, then the field $ L = K ((Y)) $ of formal series
in $ Y $ over $ K $ is infinite. We can therefore apply the
previous result which assures us that $ f $ is a product of at most 
$ M $ elements of $ D_{K ((Y)) (X)} (\f) $. But $ K ((Y)) (X) \subset
K (X) ((Y)), $ and therefore $ f \in D_{K (X) ((Y))} (\f)^M $. 
 Therefore, using lemma \ref{lelemmebonusdeDetlev}, we get that $ f \in D_{K (X)} (\f)^M $. 
$ \Box $

\subsection{Application to stably birational equivalence}

Theorem  \ref{iocentrik} allows to obtain  the
following criterion in a simple way : 

\begin{theo}
\label{zoumbalawe} Let $ \f $ and $ \psi $ be two anisotropic $ K $-forms of dimension $\geq2$,
then the following assertions are equivalent : 
\begin{enumerate}
\item $ \f \stb \psi $ ; 

\item for 
every field extension $F/K$, $ D_{F} (\psi) D_{F} (\psi) = D_{F} (\f) D_{F} (\f) $ ; 

\item for 
every field extension $F/K$, $ N_{F} (\psi) = N_{F} (\f) $.
\end{enumerate}

\noindent If  $ \f $ and $ \psi $ represent 1, these assertions
are also equivalent to
\begin{enumerate}
\item [(4)] for 
every field extension $F/K$, $ T_{F} (\psi) = T_{F} (\f) $. \end{enumerate}
\end{theo} 

\section{Link between the isotropy of $\f$ over $ K (\psi) $
and the isotropy of $ \f \otimes \pi $ over $ K(\psi \otimes \pi) $ }

A Pfister multiple is given by a pair $ (\f,\pi) $ where
$\pi\in P_{n} K $ and $\f $ is a $ K $-form. Given a Pfister form
$\pi\in P_{n} K $, we can ask if the isotropy of
$\f $ over $ K (\psi) $ has an influence on the behavior of
$\f\otimes\pi $ over $ K (\psi\otimes\pi) $. The same question arises
also for the stably birational equivalence and
hyperbolicity :  are they preserved by multiplying by
a Pfister form ? The case of hyperbolicity has been treated in
\cite{unbelareticledefotzegraldquejeneconnaisaispas} by 
Fitzgerald who proved the following result (\cite{unbelareticledefotzegraldquejeneconnaisaispas}, Theorem
3.2.) :  

\begin{theo}\label{ilestfortcefitzgerald}
Let $\f $ and $\psi $ be two $ K $-forms of dimensions $\geq 2 $, if
$\f $ is hyperbolic over $ K (\psi) $, then for every form
$\pi\in GP_{n} K $, $\pi\otimes\f $ is hyperbolic over
$ K (\pi\otimes\psi) $. 
\end{theo}

We propose to prove the following result, which was to
our knowledge not written down explicitly, although it was
already known as a consequence of a famous result by Elman and
Lam (\cite{u1}, proposition 2.2). Note that a stronger version appears in \cite{james}, corollary 2.9. 

\begin{theo}
\label{kih}Let $\f $ and $\psi $ be two $ K $-forms of dimensions
$\geq 2 $, if $\f $ is isotropic over $ K (\psi) $, then for any
form $\pi\in GP_{n} K $, $\pi\otimes\f $ is isotropic over
$ K (\pi\otimes\psi) $. 
\end{theo}

\noindent The proof will need the following result.

\begin{lem}\label{doublinclu}  Let $\f $ and $\psi $ be two anisotropic $ K $-forms
 of dimensions $\geq 2 $, then the following  assertions
are equivalent : 
\begin{enumerate}
\item $\f $ is isotropic over $ K (\psi) $ ; \item for 
every field extension $F/K$, $ D_{F} (\psi) D_{F} (\psi)\subset D_{F} (\f) D_{F} (\f) $ ; 
\item $\f\perp T\f $ is isotropic over $ K (T) (\psi\perp T\psi) $ ; 
\item $\f\perp T\f $ is isotropic over $ K ((T)) (\psi\perp T\psi) $ ; 
\item for every $ n\geq 0 $, $\l\l
X_{1},\ldots, X_{n}\r\r\otimes\f_{K (X_{1},\ldots, X_{n}) (\l\l
X_{1},\ldots, X_{n}\r\r\otimes\psi)} $ is isotropic ; \item for
every $ n\geq 0 $, $\l\l
X_{1},\ldots, X_{n}\r\r\otimes\f_{K ((X_{1}))\ldots ((X_{n})) (\l\l
X_{1},\ldots, X_{n}\r\r\otimes\psi)} $
is isotropic.\\
\end{enumerate} 
\end{lem}

\noindent\textbf{Proof  :  }$ (1)\Rightarrow(2) $ follows from  lemma 
 \ref{pouet} and $ (2)\Leftrightarrow(3)\Leftrightarrow(4) $ follows from \cite{izhboldin}, lemma  5.3. $ (4)\Rightarrow(6) $ and $ (3)\Rightarrow(5) $
appear by iterations of the implications $ (1)\Rightarrow(3) $ and
$ (1)\Rightarrow(4) $, noting that $\l\l-T\r\r\otimes\f $ becomes
isotropic over $ K (T) (\l\l-T\r\r\otimes\psi) $ in $ (3) $, and over
$ K ((T)) (\l\l-T\r\r\otimes\psi) $ in $ (4) $, but as $ K (T) = K (-T) $
and $ K ((T)) = K ((-T)) $,  we get 
$ (4)\Rightarrow(6) $ and $ (3)\Rightarrow(5) $. $ (4) $ and $ (3) $ are special cases
of $ (5) $ and $ (6) $. So we have $ (1)\Rightarrow(2) $ and
$ (2)\Leftrightarrow(3)\Leftrightarrow(4)\Leftrightarrow(5)\Leftrightarrow(6) $.
To conclude, we can use theorem \ref{iocentrik} to get $ (2)\Rightarrow(1) $. $\Box$\\

\noindent\textbf{Proof of theorem   \ref{kih} :  } We may assume that $ \pi \in P_{n} K $ and that
the form $ \pi \otimes \f $ is anisotropic over $ K $. So let us write
$ \pi \simeq \l \l a_{1}, \ldots, a_{n} \r \r $. As $ \f $ is isotropic
on $ K (\psi) $, the form $ \l \l X_{1}, \ldots, X_{n} \r \r \otimes \f $ is
isotropic over $ K (X_{1}, \ldots, X_{n}) (\l \l
X_{1}, \ldots, X_{n} \r \r \otimes \psi) $ according to  lemma 
 \ref{doublinclu}.\\
Let us write $ T = X_n-a_{n} $, then
$ K (X_{1}, \ldots, X_n) = K (X_{1}, \ldots, X_{n-1}, T) $. Consider the
completion $ F $ of $ K (X_{1}, \ldots, X_{n-1}, T $) for the $ T $-adic valuation.\\
Then $ F = K (X_{1}, \ldots, X_{n-1}) ((T)) $ and its residue class field
is $ f = K (X_{1}, \ldots, X_{n-1}) $. Then $ X_n \equiv a_n \mod
 T $ and the forms $ \pi \otimes \f $ and $ \pi \otimes \psi $ are
unimodular with residual forms over $ f $ given by
$ \l \l X_{1}, \ldots, X_{n-1}, a_n \r \r \otimes \f $ (anisotropic by
hypothesis) and $ \l \l X_{1}, \ldots, X_{n-1}, a_n \r \r \otimes \psi $. We
can then apply \cite{izhboldin}, lemma  5.3 and deduce that
$ \l \l X_{1}, \ldots, X_{n-1}, a_{n} \r \r \otimes \f $ is isotropic over
$ K (X_{1}, \ldots, X_{n-1}) (\l \l
X_{1}, \ldots, X_{n-1}, a_{n} \r \r \otimes \psi) $. By repeating this process, we conclude
that $ \pi \otimes \f $ is
isotropic over $ K (\pi \otimes \psi) $. 
$\Box$\\

Note that, more precisely, we have  shown the following result  :  

\begin{theo}\label{decadix}Let $ \f $ and $ \psi $ be two anisotropic $ K $-forms
of dimensions $ \geq 2 $. $ \f $ is
isotropic over $ K (\psi) $ if and only if for 
every field extension $F/K$, for every $n\geq0$ and for every form $ \pi \in GP_{n} F $, $ \pi \otimes \f $
is isotropic over $ F (\pi \otimes \psi) $. 
\end{theo}

We can now deduce from the previous results the following  :

\begin{theo}
\label{ikakou} 
Let $ \f $ and $ \psi $ be two anisotropic $ K $-forms
of dimensions $ \geq 2 $, then the following assertions are
equivalent : 
\begin{enumerate}
\item $ \f \stb \psi $ over $ K $ ; 

\item for every field extension $F/K$,
$ D_{F} (\psi) D_{F} (\psi) = D_{F} (\f) D_{F} (\f) $ ;

  \item for every
$ n \geq1,$  $$\l \l X_{1}, \ldots, X_{n} \r \r \otimes \f \stb \l \l
X_{1}, \ldots, X_{n} \r \r \otimes \psi $$
over $ K (X_{1}, \ldots, X_{n})$ ; 

\item for every field extension $F/K$, for every $ n \geq 1 $ and for
every $ \tau \in GP_{n} (F) $, $$ \f \otimes \tau \stb \psi \otimes \tau. $$
\end{enumerate}

\end{theo}

\begin{corol} \label{giukfdbhflkgbhgblhikblitrbh}
Let $ \f $ and $ \psi $ be two anisotropic $ K $-forms of dimensions $ \geq 2 $
such that $ \f \stb \psi $, then for every $ n \geq 1 $ and for every
form $ \pi \in GP_{n} K $, $ \pi \otimes \f \stb \pi \otimes \psi $.
\end{corol}

\section{stable birational equivalence and Clifford groups}

\subsection{Special Clifford group and spinor norm}

Let $ \f $ be a quadratic form over $ K $, $ C (\f) $ its Clifford algebra, $ C_{0} (\f) $ the even part of its Clifford algebra
and $ C_{1} (\f) $ its odd part. We borrow the following notations and results from \cite{scharlau},
chapter 1 and 9. Considering $ \f :  V \rightarrow K $,
we identify $ V $ with the subspace $ \overline{V} $ of $ C (\f) $ given
by the injection
$ \overline{i} :  \left \lbrace \begin{array}{lll}
V & \rightarrow & C (\f) \\
x & \mapsto & \overline{i} (x) = \overline{x} \end{array} \right .. $\\
The group $ O_{K} (\f) = \lbrace isometries \
\sigma :  (V, \f) \rightarrow (V, \f) \rbrace $ is called
the orthogonal group  and
$ SO_{K} (\f) = \lbrace \sigma \in O_{K} (\f) |  \det \sigma = 1 \rbrace $ the
special orthogonal group.  Any $ \sigma \in O_{K} (\f) $ defines an algebra
automorphism  $ c (\sigma) $ of $ C (\f) $. This induces a
group homomorphism 
$c :  O_{K} (\f) \rightarrow \mathrm{Aut} (C (\f)). $ Consider
the automorphism $ \gamma = c (-\mathrm{id}) $. It defines a
grading of $ C (\f) $ by $ \gamma _{\vert C_{0}} = id $ and
$ \gamma _{\vert C_{1}} = -id $.\\ The Clifford group \label{clifgroups} and the special Clifford group of $ \f $ over $ K $  are respectively defined by 
$ \Gamma_{K} (\f) = \lbrace \alpha \in C (\f) | \alpha \ invertible \ and \ \gamma (\alpha) V \alpha^{-1} = V \in \ C(\f) \rbrace $ and 
$ S \Gamma_{K} (\f) = \lbrace \alpha \in C_{0} (\f) | \alpha \ 
invertible \ and \ \alpha V \alpha^{-1} = V \in \ C(\f) \rbrace. $\\ 
We then get the following commutative diagram (\cite{spinor},
\cite{bourbakiquad} 9, $ n \raisebox{1ex}{\scriptsize o} 4 $,
theorems 2 and 3, $ n \raisebox{1ex}{\scriptsize o} 5 $ theorem 4) : 
\begin{equation*}
\begin{CD}
 K^{\ast} @ >>> S \Gamma_{K} (\f) @ >>> SO_{K} (\f) @ >>> 1 \\
  && @VV \sn V @VV \sn V \\
 && K^{\ast} @ >>> K^{\ast} / ({K^{\ast}}^{2}) @ >>> 1
\end{CD}
\end{equation*}
where $ \sn $ is the spinor norm,
given by
$ \sn (\alpha) = \alpha^{t} \alpha $ and $ S \Gamma_{K} (\f) $ acts on $ V $
by $ \alpha (v) = \alpha v \alpha^{-1} $. %Note that $ \sn $ comes from a multiplicative group homomorphism 
% $ N :  \Gamma_{K} (\f) \rightarrow K^{\ast} $.  
We can then define the following group : 
${S \Gamma_{K} (\f)}^{\diamondsuit} = S \Gamma_{K} (\f) / \ker (\sn). $
%\begin{center} and \end{center}
%$$ SO_{K} (\f)^{\diamondsuit} = SO_{K} (\f) / \ker (\sn). $$
We will denote by $ \sn^{\diamondsuit} : 
{S \Gamma_{K} (\f)}^{\diamondsuit} \rightarrow K^{\ast} $
the  injective homomorphism induced by $ \sn $. The image of 
$ \alpha \in S \Gamma_{K} (\f) $ in
${S \Gamma_{K} (\f)}^{\diamondsuit} $ under the canonical surjection
will be denoted by $ \alpha^{\diamondsuit} $.

\subsection{Structure of ${S \Gamma_{K} (\f)}^{\diamondsuit} $ and
$ \overline{S \Gamma_{K} (\f)}^{\diamondsuit} $}

Following \cite{spinor}, we will say that an element $ s \in S \Gamma_{K} (\f) $ is plane if $ s \in S \Gamma_{K} (\f_{0}) $ for $ \f_{0} \subset \f $ of dimension 2.
According to \cite{deheuvels}, Theorem IX.1,
$ S \Gamma_{K} (\f) $ is generated by its plane elements. An
element $ a \in{S \Gamma_{K} (\f)}^{\diamondsuit} $
 will be called plane if $ a = \alpha^{\diamondsuit} $ with $ \alpha $
 a plane element in $ S \Gamma_{K} (\f) $. \\
As in \cite{spinor} we will denote by
$ \overline{S \Gamma_{K} (\f)} $
 the quotient of $ S \Gamma_{K} (\f) $ by its
commutator subgroup and elements of the form
$ \alpha \beta^{-1} $
 where $ \alpha, \beta \in S \Gamma_{K} (\f) $
are plane such that $ \sn (\alpha) = \sn (\beta) $. So if we denote by
$ \overline{\sn} :  \overline{S \Gamma_{K} (\f)} \rightarrow
K^{\ast} $ the homomorphism induced by $ \sn $, we can define the group
$ \overline{S \Gamma_{K} (\f)}^{\diamondsuit} =
\overline{S \Gamma_{K} (\f)} / \ker (\overline{\sn}) $ and the 
 induced homomorphism $ \overline{\sn
}^{\diamondsuit} : 
\overline{S \Gamma_{K} (\f)}^{\diamondsuit} \rightarrow
K^{\ast} $ . The image of an element $ \alpha $ of
${S \Gamma_{K} (\f)}$ in
$\overline{S \Gamma_{K} (\f)}^{\diamondsuit}$ by the canonical surjection
will be denoted by $ \overline{\alpha}^{\diamondsuit} $. Finally, an element
$ a \in \overline{S \Gamma_{K} (\f)}^{\diamondsuit} $ will be called plane
if $ a = \overline{\alpha}^{\diamondsuit} $ with $ \alpha $ plane in $ S \Gamma_{K} (\f) $. 
Then, let us write $ P_{K} (\f) = \lbrace{\sn (\alpha) \ | \ \alpha \in
S \Gamma_{K} (\f) \ plane} \rbrace \subset K^{\ast} $ and consider
$ \widetilde{N} _{K} (\f) $, the multiplicative group 
generated by $ P_{K} (\f) $ in $ K^{\ast} $. Under these conditions, the
structures of ${S \Gamma_{K} (\f)}^{\diamondsuit} $ and
$ \overline{S \Gamma_{K} (\f)}^{\diamondsuit} $ %and $ SO_{K} (\f)^{\diamondsuit} $ 
are given by the following statement :

\begin{pro} \hfill
\label{icoptere}
\begin{enumerate}
\item ${S \Gamma_{K} (\f)}^{\diamondsuit} \simeq\widetilde{N} _{K} (\f) $.

\item $ \overline{S \Gamma_{K} (\f)}^{\diamondsuit} \simeq\widetilde{N} _{K} (\f) $.

%\item $ SO_{K} (\f)^{\diamondsuit} \simeq\widetilde{N} _{K} (\f) /{K^{\ast}}^{2} $.
\end{enumerate}
\end{pro}

\noindent \textbf{Proof  : }\\
(1) Consider
$ \Theta _{\f} :  \widetilde{N} _{K} (\f) \rightarrow{S \Gamma_{K} (\f)}^{\diamondsuit} $
 defined by
$ \Theta _{\f} (\Pi_{i} \sn (\alpha_{i})) = \Pi_{i} \alpha_{i}^{\diamondsuit
}. $ $ \Theta _{\f} $ is then well defined on
$ \widetilde{N} _{K} (\f) $ by definition of
${S \Gamma_{K} (\f)}^{\diamondsuit} $. Indeed, if
$ \Pi_{i} \sn (\alpha_{i}) = \Pi_{k} \sn (\beta_{k}) $, then
$ \Pi_{i} \sn^{\diamondsuit} (\alpha_{i}^{\diamondsuit}) =
\Pi_{k} \sn^{\diamondsuit} (\beta_{k}^{\diamondsuit}) $ by
definition of $ \sn^{\diamondsuit} $, and therefore, as
$ \sn^{\diamondsuit} $ is injective,
$ \Pi_{i} \alpha_{i}^{\diamondsuit} = \Pi_{k} \beta_{k}^{\diamondsuit
} $.
 This means exactly that
$ \Theta _{\f} (\Pi_{i} \sn (\alpha_{i})) = \Theta _{\f} (\Pi_{k} \sn (\beta_{k})). $\\
The map $ \Theta _{\f} $ is surjective since
$ S \Gamma_{K} (\f) $, and therefore ${S \Gamma_{K} (\f)}^{\diamondsuit} $,
is generated by their plane elements. This assures us
that, for every $ x \in{S \Gamma_{K} (\f)}^{\diamondsuit} $, we can
write
$ x = \Pi_{i}{\alpha_{i}}^{\diamondsuit} =
\Theta_{\f} (\Pi_{i} \sn (\alpha_{i})) $ where $ \Pi_{i} \sn (\alpha_{i}) \in
\widetilde{N} _{F} (\f) $ and the $\alpha_{i}$ are plane.\\
$ \Theta _{\f} $  is also injective. Indeed, suppose
$ \Theta _{\f} (x) = 1^{\diamondsuit} $,
 and write
$ x = \Pi_{i} \sn (\alpha_{i}). $ Consequently, we will have
 $ \Theta _{\f} (\Pi_{i} \sn (\alpha_{i})) = 1^{\diamondsuit} $ and therefore,
by definition of $ \Theta _{\f} $,
$ \Pi_{i} \alpha_{i}^{\diamondsuit} = 1^{\diamondsuit} $. So  $ \sn^{\diamondsuit} (\Pi_{i}{\alpha_{i}}^{\diamondsuit}) =
\sn^{\diamondsuit} (1^{\diamondsuit}) = $ 1 but then
$ \Pi_{i} \sn (\alpha_{i}) = \Pi_{i} \sn^{\diamondsuit
} ({\alpha_{i}}^{\diamondsuit}) = 1, $
 which gives
$ x = $ 1. So the map $ \Theta _{\f} $ is a group isomorphism.\\
(2) This part is shown as the previous one, defining
$ \overline{\Theta _{\f}} :  \widetilde{N} _{K} (\f)
\rightarrow \overline{S \Gamma_{K} (\f)}^{\diamondsuit}$ by $\overline{\Theta _{\f}} (\Pi_{i} \sn (\alpha_{i})) = \Pi_{i} \overline{\alpha_{i}}^{\diamondsuit
}. $
 Injectivity
of $ \overline{\sn}^{\diamondsuit} $ will once again assure us that
the homomorphism is well defined and the rest is similar.$ \Box $ 
%3. \cite{algebraicspinor} page 53, II.3.6 and lemma 8 (i) from \cite{spinor} show that $ D_{K} (\f) D_{K} (\f) = P_{K} (\f) $. 

\subsection{Stable birational equivalence and group isomorphisms}

Rost showed in lemma  8 of \cite{spinor} that
for every $ K $-quadratic form $ \f $,
$ D_{K} (\f) D_{K} (\f) = P_{K} (\f) $ where $ P_{K} (\f) $ is the set
previously defined. As a consequence, the sets
$ \widetilde{N} _{K} (\f) $ and $ N_{K} (\f) $ coincide. 
Then, we can use theorem  \ref{zoumbalawe}  to 
get the following criterion for stable birational equivalence of quadratic forms :  

\begin{theo} \label{letheoremsuivantrost}
Let $ \f $, $ \psi $ be two anisotropic $ K $-quadratic forms, then the following assertions are equivalent : 
\begin{enumerate}
\item $ \f \stb \psi $ ; 

\item for every field extension $F/K$, there exists an group isomorphism
$$ \xi_{F} :  S{\Gamma_{F} (\f)}^{\diamondsuit} \rightarrow
S{\Gamma_{F} (\psi)}^{\diamondsuit} $$ such that for every $ x \in
S{\Gamma_{F} (\f)}^{\diamondsuit}, $ $ \sn^{\diamondsuit} (x) =
\sn^{\diamondsuit} (\xi_{F} (x)) \ ;  \ $

\item for every field extension $F/K$, there exists a group isomorphism
$$ \overline{\xi} _{F} : 
 \overline{S \Gamma_{F} (\f)}^{\diamondsuit} \rightarrow
 \overline{S \Gamma_{F} (\psi)}^{\diamondsuit} $$
 such that for every $ \overline{x} \in
\overline{S \Gamma_{F} (\f)}^{\diamondsuit}, $ $ \overline{\sn
}^{\diamondsuit} (\overline{x}) = \overline{\sn}^{\diamondsuit
} (\overline{\xi} _{F} (\overline{x})). $
\end{enumerate}
\end{theo} 

\noindent\textbf{Proof : }
$(1)\Rightarrow(2)$ Theorem  \ref{zoumbalawe} shows that $ \f \stb \psi $ if and
only if $ N_{F} (\f) = N_{F} (\psi) $ for 
every field extension $F/K$ and therefore if and only if
$ \widetilde{N} _{F} (\f) = \widetilde{N} _{F} (\psi) $ for 
every field extension $F/K$. Then, since Proposition \ref{icoptere} (1) gives 
${S \Gamma_{K} (\f)}^{\diamondsuit} \simeq\widetilde{N} _{K} (\f) $
and ${S \Gamma_{K} (\psi)}^{\diamondsuit} \simeq\widetilde{N} _{K} (\psi) $, we have
${S \Gamma_{K} (\f)}^{\diamondsuit} \simeq {S \Gamma_{K} (\psi)}^{\diamondsuit}$.
More precisely,   the group
isomorphism constructed in the proof of the proposition  \ref{icoptere} (1) actually satisfies the additional following property
 : 
$ \forall a \in{S \Gamma_{F} (\f)}^{\diamondsuit}, $
$ \sn^{\diamondsuit} (a) = \sn^{\diamondsuit} (\xi (a)).$\\
$(2)\Rightarrow(1)$  Consider $\Theta_\f$ as defined in the proof of \ref{icoptere}. Since the isomorphism between
${S \Gamma_{F} (\f)}^{\diamondsuit} $ and
$ S{\Gamma_{F} (\psi)}^{\diamondsuit} $ satisfies $ \forall a \in{S \Gamma_{F} (\f)}^{\diamondsuit}, $
$ \sn^{\diamondsuit} (a) = \sn^{\diamondsuit} (\xi (a))$, then,  for every $a \in N_{F} (\f),$ we have
$ \sn^{\diamondsuit} (\xi (\Theta _{\f} (a))) = a \in N_{F} (\psi).$ As a consequence, we have $ N_{F} (\f) = N_{F} (\psi) $.
Therefore, using Theorem \ref{zoumbalawe} we get $ \f \stb \psi $.\\
$(1)\Leftrightarrow(3)$ is similar as $(1)\Leftrightarrow(2)$, replacing ${S \Gamma_{F} (\f)}^{\diamondsuit}$ by ${\overline{S \Gamma_{F} (\f)}}^{\diamondsuit}$,  $\sn^{\diamondsuit}$ by $\overline{\sn}^{\diamondsuit}$, $\Theta _{\f}$ by $\overline{\Theta _{\f}}$ and 
$\xi _{F}$ by $\overline{\xi} _{F}$.
$\Box$

%\subsection*{acknowledgements}

% ----------------------------------------------------------------
\bibliographystyle{amsplain}
\bibliography{critstabireqbiblio}
\end{document}